\documentclass[12pt]{article}
\usepackage{amsthm,amsfonts,amssymb,amscd,amsmath}

\begin{document}

\begin{center}
\Large \bf On the global log canonical threshold\\
of Fano complete intersections
\end{center}
\vspace{1cm}

\centerline{Thomas Eckl and Aleksandr Pukhlikov}

\parshape=1
3cm 10cm \noindent {\small \quad \quad \quad
\quad\quad\quad\quad\quad\quad\quad {\bf }\newline We show that
the global log canonical threshold of generic Fano complete
intersections of index 1 and codimension $k$ in ${\mathbb
P}^{M+k}$ is equal to 1 if $M\geqslant 3k+4$ and the highest
degree of defining equations is at least 8. This improves the
earlier result where the inequality $M\geqslant 4k+1$ was
required, so the class of Fano complete intersections covered by
our theorem is considerably larger. The theorem implies, in
particular, that the Fano complete intersections satisfying our
assumptions admit a K\" ahler-Einstein metric. We also show the
existence of K\" ahler-Einstein metrics for a new finite set of
families of Fano complete intersections.

Bibliography: 9 titles.} \vspace{1cm}

\vspace{1cm}

14E05, 14E07, 14J45\vspace{1cm}

\noindent {\bf 1. The canonical and log canonical thresholds.}
Consider a smooth Fano variety $X$, such that $\mathop{\rm
Pic}X={\mathbb Z}K_X$, $D\sim-nK_X$ an effective divisor, $n\geq
1$.\vspace{0.1cm}

{\bf Definition 1.} The pair $(X,\alpha D)$, where
$\alpha\in{\mathbb Q}$, is {\it not canonical} (respectively, {\it
not log canonical}), if there exists a prime divisor $E$ over $X$
such that the inequality
$$
\alpha\mathop{\rm ord}\nolimits_E D>a(E)
$$
(respectively, $\alpha\mathop{\rm ord}_E D>a(E)+1$) is satisfied,
where $a(E)=a(X,E)$ is the discrepancy of $E$ with respect to the
model $X$.\vspace{0.1cm}

Explicitly, this means that there are a birational morphism
$\varphi\colon\widetilde X\to X$, where $\widetilde X$ is a smooth
projective variety, and a prime $\varphi$-exceptional divisor
$E\subset\widetilde X$ such that
$$
\alpha\mathop{\rm ord}\nolimits_E\varphi^*D>a(E)
$$
(respectively, $\alpha\mathop{\rm
ord}_E\varphi^*D>a(E)+1$).\vspace{0.1cm}

{\bf Definition 2.} The global canonical (respectively, log
canonical) threshold of the variety $X$ is defined by the equality
$$
\mathop{\rm ct}(X)=\mathop{\rm sup}\{\lambda\in{\mathbb
Q}_+\,|\,\mbox{the pair}\,\, (X,\frac{\lambda}{n}D)\,\,\mbox{is
canonical for all}\,\,D\in |-nK_X|\},
$$
(respectively,
$$
\mathop{\rm lct}(X)=\mathop{\rm sup}\{\lambda\in{\mathbb
Q}_+\,|\,\mbox{the pair}\,\,(X,\frac{\lambda}{n}D)\,\,\mbox{is log
canonical for all}\,\,D\in |-nK_X|\}).
$$
(Note that as $D\sim-nK_X$, the integer $n\geqslant 1$ depends on
the effective divisor $D$.)\vspace{0.1cm}

The importance of canonical and log canonical thresholds comes
from their applications to complex differential geometry and to
certain problems of higher-dimensional birational geometry.

In \cite{DeKo2001,Nadel90,Tian87} the following fact was
shown.\vspace{0.1cm}

{\bf Theorem 1.} {\it Assume that the inequality
$$
\mathop{\rm lct}(X)>\frac{\mathop{\rm dim}X}{\mathop{\rm dim}X+1}
$$
holds. Then on the variety $X$ there exists a K\" ahler-Einstein
metric.}\vspace{0.1cm}

The log canonical threshold is important in this
differential-geometric context because it indicates the
non-triviality of certain multiplier ideal sheaves. In their
analytic interpretation, these multiplier ideal sheaves in turn
measure the failure of a priori estimates sufficient to solve a
Monge-Amp\`ere equation for a K\"ahler-Einstein metric.

{\bf Definition 3.} The {\it mobile canonical threshold} of a
variety $X$, which is denoted by the symbol $\mathop{\rm mct}(X)$,
is the supremum of such $\lambda\in {\mathbb Q}_+$ that the pair
$(X,\frac{\lambda}{n}D)$ is canonical for a generic divisor $D\in
\Sigma$ of any mobile linear system $\Sigma\subset |-nK_X|$ (that
is to say, any system $\Sigma$ with no fixed
components).\vspace{0.1cm}

In \cite{Pukh05} the following fact was shown.\vspace{0.1cm}

{\bf Theorem 2.} {\it Assume that primitive Fano varieties
$F_1,\dots,F_K$, $K\geq 2$, satisfy the conditions $\mathop{\rm
lct} (F_i)= \mathop{\rm mct}(F_i) = 1$. Then their direct
product
$$
V=F_1\times\dots\times F_K
$$
is a birationally superrigid variety. In particular,\vspace{0.1cm}

{\rm (i)} Every structure of a rationally connected fiber space on
the variety $V$ is given by a projection onto a direct factor.
More precisely, let $\beta\colon V^{\sharp}\to S^{\sharp}$ be a
rationally connected fiber space and $\chi\colon V\dashrightarrow
V^{\sharp}$ a birational map. Then there exists a subset of
indices
$$
I=\{i_1,\dots,i_k\}\subset \{1,\dots,K\}
$$
and a birational map
$$
\alpha\colon F_I=\prod\limits_{i\in I}F_i \dashrightarrow
S^{\sharp},
$$
such that the diagram
$$
\begin{array}{rcccl}
& V &\stackrel{\chi}{\dashrightarrow} & V^{\sharp}&\\
\pi_I &\downarrow & &\downarrow &\beta\\
& F_I & \stackrel{\alpha}{\dashrightarrow}& S^{\sharp}&
\end{array}
$$
commutes, that is, $\beta\circ\chi=\alpha\circ\pi_I$, where
$$
\pi_I\colon\prod\limits^K_{i=1}F_i\to \prod\limits_{i\in I}F_i
$$
is the natural projection onto a direct factor.\vspace{0.1cm}

{\rm (ii)} Let $V^{\sharp}$ be a variety with ${\mathbb
Q}$-factorial terminal singularities, satisfying the condition
$$
\mathop{\rm dim}\nolimits_{\mathbb Q}(\mathop{\rm
Pic}V^{\sharp}\otimes{\mathbb Q})\leqslant K,
$$
and $\chi\colon V\dashrightarrow V^{\sharp}$ a birational map.
Then $\chi$ is a (biregular) isomorphism.\vspace{0.1cm}

{\rm (iii)} The groups of birational and biregular self-maps of
the variety $V$ coincide:
$$
\mathop{\rm Bir}V=\mathop{\rm Aut}V.
$$
In particular, the group $\mathop{\rm Bir}V$ is
finite.\vspace{0.1cm}

{\rm (iv)} The variety $V$ admits no structures of a fibration
into rationally connected varieties of dimension strictly smaller
than $\mathop{\rm min}\{\mathop{\rm dim}F_i\}$. In particular, $V$
admits no structures of a conic bundle or a fibration into
rational surfaces.\vspace{0.1cm}

{\rm (v)} The variety $V$ is non-rational.}\vspace{0.1cm}

For the precise definition of birational (super)rigidity, a
discussion of its properties and examples of birationally
(super)rigid varieties, see \cite{Pukh13a}.\vspace{0.3cm}


{\bf 2. Fano complete intersections.} Fix an integer $k\geqslant
2$. Consider an arbitrary system $(d_1,\dots,d_k)$ of positive
integers, satisfying the condition $d_k\geqslant\dots\geqslant
d_1\geqslant 2$. Set $M=d_1+\dots+d_k-k$. Fix the complex
projective space ${\mathbb P}={\mathbb P}^{M+k}$ and consider the
family ${\cal F}(d_1,\dots,d_k)$ of non-singular complete
intersections $V$ of the type $d_1\cdot\dots\cdot d_k$ in
${\mathbb P}$:
$$
V=F_1\cap\dots\cap F_k\subset{\mathbb P},
$$
Here $F_i\subset{\mathbb P}$ is a hypersurface of degree $d_i$,
and $\mathop{\rm codim}V=k$.\vspace{0.1cm}

The following two theorems collect the known information about the
global canonical and log canonical thresholds of Fano complete
intersections as above.\vspace{0.1cm}

{\bf Theorem 3.} {\it Assume that $M\geqslant 4k+1$ and
$d_k\geqslant 8$. Then for a generic (in the sense of Zariski
topology on the space ${\cal F}(d_1,\dots,d_k)$) variety $V\in
{\cal F}(d_1,\dots,d_k)$ the equality $\mathop{\rm ct} (V)=1$
holds.}\vspace{0.1cm}

{\bf Proof:} see \cite[Section 3]{Pukh06b}.\vspace{0.1cm}

Thus under the assumptions of Theorem 3 on $V$ exists a K\"
ahler-Einstein metric. Besides, since $\mathop{\rm lct}
(V)\geqslant \mathop{\rm ct} (V)$ and $\mathop{\rm mct}
(V)\geqslant \mathop{\rm ct} (V)$, the variety $V$ satisfies the
assumptions of Theorem 2 and for that reason can be used as a
direct factor in birationally superrigid Fano direct
products.\vspace{0.1cm}

{\bf Theorem 4.} {\it Assume that $M\geqslant 4k+1$ and any of the
following conditions holds:

{\rm (i)} $d_k=d_{k-1}=7$ and $M\leq 47$,

{\rm (ii)}  $d_k=7$, $d_{k-1}\leq 6$ and $M\leq 19$.

{\rm (iii)} $k=2$, $d_1=d_2=6$, $M=10$.

\noindent Then the canonical threshold $\mathop{\rm ct}(V)$ (and
hence also the log canonical threshold $\mathop{\rm lct}(V)$) of a
generic variety $V\in{\cal F}(d_1,\dots,d_k)$ satisfies the
inequality}
$$
\mathop{\rm ct}(V)>\frac{M}{M+1}.
$$

{\bf Proof:} see \cite{Pukh10c}.\vspace{0.1cm}

Therefore,  on a general variety $V\in {\cal F}(d_1,\dots,d_k)$,
satisfying one of the conditions listed in Theorem 4 there exists
a K\"ahler-Einstein metric.\vspace{0.1cm}

The aim of the present note is to improve the claims of Theorems 3
and 4, extending them to a larger class of Fano complete
intersections of index 1. We will show the following two
facts.\vspace{0.1cm}

{\bf Theorem 5.} {\it Assume that $M\geqslant 3k+4$ and
$d_k\geqslant 8$. Then for a generic (in the sense of Zariski
topology on the space ${\cal F}(d_1,\dots,d_k)$) variety $V\in
{\cal F}(d_1,\dots,d_k)$ the equality $\mathop{\rm ct} (V)=1$
holds.}\vspace{0.1cm}

{\bf Theorem 6.} {\it Assume that $M\geqslant 3k+4$ and any of the
two following conditions holds:

{\rm (i)} $d_k=d_{k-1}=7$ and $M\leqslant 47$,

{\rm (ii)}  $d_k=7$, $d_{k-1}\leqslant 6$ and $M\leqslant 19$.

{\rm (iii)} $k=2$, $d_1=d_2=6$, $M=10$.

\noindent Then the canonical threshold $\mathop{\rm ct}(V)$  (and
hence also the log canonical threshold $\mathop{\rm lct}(V)$) of a
generic variety $V\in{\cal F}(d_1,\dots,d_k)$ satisfies the
inequality}
$$
\mathop{\rm ct}(V)>\frac{M}{M+1}.
$$

{\bf Remark 1.} Theorem 5 covers a considerably larger class of
Fano complete intersections than Theorem 3. The same is true for
the part (i) of Theorems 6 and 4. For the part (iii), Theorem 6
gives nothing new compared to Theorem 4. As for the part (ii),
Theorem 6 gives the existence of the K\" ahler-Einstein metric for
the following 7 families of Fano complete intersections that do
not fit into the assumptions of Theorem 4, all of them in
${\mathbb P}^{24}$:
$$
\begin{array}{c}
{\cal F}(2,5,5,5,7),\quad {\cal F}(2,4,5,6,7),\quad {\cal
F}(2,3,6,6,7),\quad {\cal
F}(3,3,5,6,7), \\
{\cal F}(3,4,5,5,7),\quad {\cal F}(3,4,4,6,7),\quad {\cal
F}(4,4,4,5,7).
\end{array}
$$
\vspace{0.3cm}


{\bf 3. The conditions of general position.} Now we will explain
what we mean by the genericity of a Fano complete intersection
$V\in {\cal F}(d_1,\dots,d_k)$. Fix any point $o\in V$ and let
$(z_1,\dots,z_{M+k})$ be a system of affine coordinates on
${\mathbb P}$ with the origin at the point $o$,
$$
f_i=q_{i,1}+\dots+q_{i,d_i}
$$
the equation of the hypersurface $F_i$ with respect to
(non-homogeneous) coordinates $z_*$, decomposed into homogeneous
in $z_*$ components $q_{i,j}$, $\mathop{\rm deg}q_{i,j}=j$. Since
$V$ is a non-singular variety, the system of linear equations
$$
q_{1,1}=\dots=q_{k,1}=0
$$
defines a linear subspace $T_oV\subset{\mathbb C}^{M+k}$ of
codimension $k$, the tangent space to the variety $V$ at the point
$o$. We define a finite set of pairs $I\subset{\mathbb
Z}_+\times{\mathbb Z}_+$ in the following way: if $d_{k-1}=d_k$,
then we set
$$
I=\{(i,j)\,|\,1\leqslant i\leq k,1\leqslant j\leqslant
d_i,(i,j)\not\in\{(k,d_k),(k-1,d_{k-1})\}\},
$$
and if $d_{k-1}\leq d_k-1$, then we set
$$
I=\{(i,j)\,|\,1\leqslant i\leqslant k,1\leqslant j\leqslant
d_i,(i,j)\not\in\{(k,d_k),(k,d_k-1)\}\}.
$$

{\bf Definition 4.} We say that the complete intersection $V$ is
{\it regular at the point} $o$, if for any linear form $l(z_*)$,
not vanishing identically on the subspace $T_oV$, the set of
homogeneous polynomials
$$
\{l\}\cup\{q_{i,j}\,|\,(i,j)\in I\}
$$
forms a regular sequence in ${\cal O}_{o,{\mathbb P}}$, that is,
the system of equations in ${\mathbb C}^{M+k}$
$$
l=0,\quad q_{i,j}=0,(i,j)\in I,
$$
defines a subset of codimension $\# I+1$. Finally, we say that
the complete intersection $V$ is {\it regular}, if it is regular
at every point.\vspace{0.1cm}

In Sec. 4 below we show the following fact.\vspace{0.1cm}

{\bf Theorem 7.} {\it For $M\geqslant 3k+4$ there exists a
non-empty Zariski open subset
$$
{\cal F}_{\rm reg}(d_1,\dots,d_k)\subset {\cal F}(d_1,\dots,d_k),
$$
such that any variety $V\in {\cal F}_{\rm reg}(d_1,\dots,d_k)$ is
regular.}\vspace{0.1cm}

By genericity in Theorems 3-6 we mean the regularity. For that
reason, the main results of this paper (Theorems 5,6) are
essentially dependent on Theorem 7 which allows us to assume
regularity of the complete intersection $V$.\vspace{0.3cm}


{\bf 4. Proof of the regularity conditions.} Let us show Theorem
7. The proof proceeds in a series of reduction steps and case
distinctions, through the estimates (\ref{incident-est-eq}) -
(\ref{k-M-3-est-eq}).\vspace{0.1cm}

Let the incidence variety $\mathcal{V} \subset
H^0(\mathbb{P}^{M+k}, \mathcal{O}_{\mathbb{P}^{M+k}}(1)) \times
\mathbb{P}^{M+k} \times \mathcal{F}(d_1, \ldots, d_k)$ consist of
tuples $(L, o, F_1, \ldots, F_k)$ such that $o \in \{L = F_1 =
\cdots = F_k = 0\}$. Let $p, q, r$ be the natural projections of
$\mathcal{V}$ to $H^0(\mathbb{P}^{M+k},
\mathcal{O}_{\mathbb{P}^{M+k}}(1))$, $\mathbb{P}^{M+k}$ resp.\
$\mathcal{F}(d_1, \ldots, d_k)$.\vspace{0.1cm}

The image of $\mathcal{V}$ under the projection $p \times q$ is
the incidence variety
\[ \mathcal{L} \subset H^0(\mathbb{P}^{M+k},
\mathcal{O}_{\mathbb{P}^{M+k}}(1)) \times \mathbb{P}^{M+k} \]
consisting of pairs $(L, o)$ such that $L(o) = 0$. All the $(p
\times q)$-fibers $\mathcal{V}_{L,o} \subset \mathcal{V}$ over
points $(L,o) \in \mathcal{L}$ are isomorphic to a product of
affine spaces and have codimension $k$ in $\mathcal{F}(d_1,
\ldots, d_k)$, as vanishing in $o$ imposes one linear condition on
the sections in $H^0(\mathbb{P}^N,
\mathcal{O}_{\mathbb{P}^N}(d_i))$. Hence $\mathcal{V}$ is
irreducible.\vspace{0.1cm}

Bertini's theorem shows that for a general tuple $(F_1, \ldots,
F_k) \in \mathcal{F}(d_1, \ldots, d_k)$ the algebraic subset
$\{F_1 = \cdots = F_k = 0\}$ is an $M$-dimensional complete
intersection. Hence the general $r$-fiber in $\mathcal{V}$ has
dimension $(M+k) + M = 2M + k$.\vspace{0.1cm}

Let $\mathcal{V}^{\mathrm{sm, nonreg}}$ be the locally
Zariski-closed subset of tuples $(L,o,F_1, \ldots, F_k) \in
\mathcal{V}$ such that $V = \{F_1 = \ldots = F_k =  0\}$ is a
smooth complete intersection of dimension $M$, but the
dehomogenisation $l$ of the linear form $L$ in affine coordinates
around $o$ and the homogeneous parts $q_{i,j}$ of the
dehomogenized $F_i$ do not satisfy the regularity condition of
Definition 4. Then the theorem is shown if the projection of the
Zariski-closure of $\mathcal{V}^{\mathrm{sm, nonreg}}$ does not
cover $\mathcal{F}(d_1, \ldots, d_k)$.\vspace{0.1cm}

To show this claim we note that every $r$-fiber of
$\mathcal{V}^{\mathrm{sm, nonreg}}$ must be at least
$k$-dimensional: If the subscheme $V = \{F_1 = \ldots = F_k =  0\}
\subset \mathbb{P}^{M+k}$ is smooth of dimension $M$ in $o \in
\mathbb{P}^{M+k}$ then we can choose homogeneous coordinates
$[Z_0: \ldots :Z_{M+k}]$ on $\mathbb{P}^{M+k}$ such that $o =
[1:0: \ldots :0]$ and $q_{1,1}= z_{M+1}, \ldots, q_{k,1} =
z_{M+k}$ in the affine coordinates $z_1, \ldots, z_{M+k}$
dehomogenized with respect to $Z_0$. Hence there is a
$k$-dimensional linear subspace of linear forms $L \in
H^0(\mathbb{P}^{M+k}, \mathcal{O}_{\mathbb{P}^{M+k}}(1))$ with the
same intersection of $\{L = 0\} \subset \mathbb{P}^{M+k}$ and the
tangent space $T_o V$, seen as a linear subspace of
$\mathbb{P}^{M+k}$. In particular, if the sequence $l, q_{1,1},
\ldots, q_{k,1}, q_{1,2}, \ldots$ is not regular for one such
linear form $L$ (dehomogenized to $l$) then the sequence will not
be regular for any other linear form in this $k$-dimensional
linear subspace. Consequently, we only need to show that
\begin{equation} \label{incident-est-eq}
\mathrm{codim}_{\mathcal{V}} \mathcal{V}^{\mathrm{sm, nonreg}}
\geq (2M+k) - k + 1 = 2M+1.
\end{equation}

The $q$-fibers $\mathcal{V}_o$ over points $o \in
\mathbb{P}^{M+k}$ are all isomorphic. Consequently it is enough to
show for all the subsets $\mathcal{V}^{\mathrm{sm, nonreg}} \cap
\mathcal{V}_o =: \mathcal{V}^{\mathrm{sm, nonreg}}_o$ locally
Zariski-closed in $\mathcal{V}_o$ that
\begin{equation}
\mathrm{codim}_{\mathcal{V}_o} \mathcal{V}^{\mathrm{sm, nonreg}}_o \geq 2M+1.
\end{equation}

Let $\mathcal{P}^N_d$ denote the vector space of homogeneous
polynomials of degree $d$ in the affine coordinates $z_1, \ldots,
z_N$, $N \leq M+k$, introduced above. Then
\[ \mathcal{V}_o \cong \mathcal{P}^{M+k}_1 \times
\prod_{i=1}^k \prod_{j=1}^{d_i} \mathcal{P}^{M+k}_j, \] by
identifying the dehomogenized sections in $H^0(\mathbb{P}^{M+k},
\mathcal{O}_{\mathbb{P}^{M+k}}(j))$ vanishing in $o$ with
$\mathcal{P}^{M+k}_j$ and decomposing the dehomogenisation $f_i =
q_{i,1} + \cdots q_{i,d_i}$ of the $F_i$ in the tuple $(L, F_1,
\ldots, F_k) \in \mathcal{V}_o$ into homogeneous parts $q_{i,j}$
of degree $j$.\vspace{0.1cm}

Consider the projection $s$ of $\mathcal{V}_o$ to
$\mathcal{P}^{M+k}_1 \times \prod_{i=1}^k (\mathcal{P}^{M+k}_1)$,
that is, to the linear form $l$ coming from $L$ and the linear
parts $q_{i,1}$ of the dehomogenized $F_i$. The map $s$ is not
defined everywhere on $\mathcal{V}_o$ but since
$\mathcal{V}^{\mathrm{sm, nonreg}}_o$ consists of tuples $(L, F_1,
\ldots, F_k)$ such that $V = \{F_1 = \ldots = F_k =  0\} \subset
\mathbb{P}^{M+k}$ is smooth of codimension $k$ in $o$ and $l_{|T_o
V} \not\equiv 0$, hence none of the $l$ and $q_{i,1}$ can be $0$,
we conclude that $s$ is defined on $\mathcal{V}^{\mathrm{sm,
nonreg}}_o$.\vspace{0.1cm}

Therefore, it is enough to show that the codimension of the
isomorphic $s$-fibers in $\mathcal{V}^{\mathrm{sm, nonreg}}_o$ is
$\geq 2M+1$. Choosing the coordinates $z_\ast$ such that $l = z_M,
q_{i,1} = z_{M+1}, \ldots, q_{k,1} = z_{M+k}$, that means to show
that the set $U(d_1, \ldots, d_k)$ of tuples of homogeneous
polynomials in variables $z_1, \ldots, z_{M-1}$ defined by
\[ \left\{ (q_{i,j})_{1 \leq i \leq k, 2 \leqslant j \leqslant d_i}\,|\,
(q_{i,j})_{(i,j) \in I, j \neq 1}\ \mathrm{is\ not\ a\ regular\
sequence\ in\  } \mathcal{O}_{0,\mathbb{A}^{M-1}} \right\}\] and
Zariski-closed in $\prod_{i=1}^k \prod_{j=2}^{d_i}
\mathcal{P}_j^{M-1}$ satisfies
\begin{equation}
\mathrm{codim}_{\prod_{i=1}^k \prod_{j=2}^{d_i}
\mathcal{P}_j^{M-1}} U(d_1, \ldots, d_k) \geqslant 2M+1.
\end{equation}

In order to obtain this estimate only the degrees of the
homogeneous polynomials $q_{i,j}$ matter. So we will discuss the
codimension of the set of tuples
\[ U := \{ (q_i)_{1 \leqslant i \leq M}\,|\, (q_i)_{1 \leqslant i
\leqslant M-2}\ \mathrm{is\ not\ a\ regular\ sequence\ in\ }
\mathcal{O}_{0,\mathbb{A}^{M-1}}  \}  \] in $\prod_{i=1}^M
\mathcal{P}_{m_i}^{M-1}$ where $2 \leqslant m_1 \leqslant \cdots
\leqslant m_M = d_k$ and we have
\[ k_d := \# \{ d_l: d_l \geqslant d, 1 \leqslant l \leqslant k\} \]
polynomials $q_i$ of degree $d$. In particular $k = k_2 \geqslant
k_3 \geqslant \cdots \geqslant k_{d_k}$ and
\[ \sum_{i=1}^M m_i = \sum_{d=2}^{d_k} k_d \cdot d =
\sum_{l=1}^k \sum_{d=2}^{d_l} d = \sum_{l=1}^k \frac{d_l(d_l+1)}{2}- k. \]

Let $Z(q_1, \ldots, q_j) \subset \mathbb{A}^{M-1}$ denote the
vanishing locus of $q_1, \ldots, q_j$ in $\mathbb{A}^{M-1}$. Since
the $q_i$ are homogeneous $Z(q_1, \ldots, q_j)$ is a cone over the
origin in $\mathbb{A}^{M-1}$, and we denote its projectivization
in $\mathbb{P}^{M-2}$ by $V(q_1, \ldots, q_j)$. In particular,
$(q_1, \ldots, q_j)$ is regular in $0 \in \mathbb{A}^{M-1}$ if and
only if
\[ \mathrm{codim}_{\mathbb{P}^{M-2}} V(q_1, \ldots, q_j) =
\mathrm{codim}_{\mathbb{A}^{M-2}} Z(q_1, \ldots, q_j) = j, \]
where the codimension is set to be the minimum of the codimensions
of the irreducible components.\vspace{0.1cm}

Consequently, $U$ is covered by locally Zariski-closed subsets
$U_j \times \prod_{i=j+1}^M \mathcal{P}_{m_i}^{M-1}$, where $1
\leqslant j \leqslant M-2$ and
\[ U_j := \{(q_1, \ldots, q_j)\,|\,
 \mathrm{codim}_{\mathbb{P}^{M-2}}V(q_1, \dots, q_j) =
 \mathrm{codim}_{\mathbb{P}^{M-2}}V(q_1, \dots, q_{j-1}) = j-1 \} \]
is a locally Zariski-closed subset of $\prod_{i=1}^j
\mathcal{P}_{m_i}^{M-1}$. Thus we need to show
\begin{equation}
\mathrm{codim}_{\prod_{i=1}^k \mathcal{P}_{m_i}^{M-1}} U = \min_{1
\leqslant j \leqslant M-2} \mathrm{codim}_{\prod_{i=1}^j
\mathcal{P}_{m_i}^{M-1}} U_j \geqslant 2M+1,
\end{equation}
and this means to verify
\begin{equation} \label{U_j-deg-2-eq}
 \mathrm{codim}_{\prod_{i=1}^j \mathcal{P}_{m_i}^{M-1}} U_j
 \geqslant 2M+1\ \mathrm{for\ } 1 \leqslant j \leqslant M-2.
\end{equation}

If $m_j=2$ (that is, $1 \leqslant j \leqslant k_2 = k$) then we
estimate the codimension of $U_j$ in $\prod_{i=1}^j
\mathcal{P}_{m_i}^{M-1}$ using the method of \cite{Pukh98b}:
Choose a general $(j-2)$-dimensional hyperplane $S \subset
\mathbb{P}^{M-2}$. Then the projection $\pi_S: \mathbb{P}^{M-2}
\rightarrow \mathbb{P}^{M-1-j}$ restricts to a finite morphism on
each of the irreducible components of $V(q_1, \ldots, q_{j-1})$
for given $q_1, \ldots, q_{j-1}$ such that
$\mathrm{codim}_{\mathbb{P}^{M-2}}V(q_1, \dots, q_{j-1}) = j-1 $.
Consequently, polynomials in $\mathcal{P}_{m_j}^{M-1} =
\mathcal{P}_2^{M-1}$ obtained as a pullback of a homogeneous
quadratic polynomial defined on $\mathbb{P}^{M-1-j}$ do not vanish
on $V(q_1, \ldots, q_{j-1})$. The linear space $W_j$ of such
pulled-back quadratic  polynomials has dimension
$\binom{M-1-j+2}{2}$. By construction, $W_j$ intersects the space
of polynomials in $\mathcal{P}_2^{M-1}$ vanishing on one of the
irreducible components of $V(q_1, \ldots, q_{j-1})$ only in $0$.
Therefore,
\begin{equation}
\mathrm{codim}_{\prod_{i=1}^j \mathcal{P}_{m_i}^{M-1}} U_j
\geqslant \binom{M-j+1}{2} \geqslant \binom{M-k+1}{2},\
    j = 1, \ldots, k.
\end{equation}
Note that this estimate also holds if $j=1$. Furthermore,
$\binom{M-k+1}{2} \geqslant 2M+1$ since $3k + 4 \leqslant M$.
Hence $(\ref{U_j-deg-2-eq})$ is shown for $j = 1, \ldots,
k$.\vspace{0.1cm}

If $m_j > 2$ (that is, $k+1 \leqslant j \leqslant M-2$) then we
estimate the codimension of $U_j$ in $\prod_{i=1}^j
\mathcal{P}_{m_i}^{M-1}$ using the method of \cite{Pukh01}: If
$(q_1, \ldots, q_j) \in U_j$ choose an irreducible component $B$
of $V(q_1, \ldots, q_j)$. By definition of $U_j$, the codimension
of $B$ in $\mathbb{P}^{M-2}$ is $j-1$. Assume that the codimension
of the linear subspace $\langle B \rangle \subset
\mathbb{P}^{M-2}$ spanned by $B$ is $b$; that means in particular
that $0 \leqslant b \leqslant j-1$.\vspace{0.1cm}

{\bf Lemma 1.} If $b < j-1$, then there are indices $1 \leqslant
i_1 < \cdots < i_{j-1-b} \leqslant j-1$ such that $B$ is an
irreducible component of
\[ V(q_{i_1}, \ldots, q_{i_{j-1-b}}) \cap \langle B \rangle. \]

{\bf Proof.} Since $\mathrm{codim}_{\mathbb{P}^{M-2}} \langle B
\rangle = b < j-1 = \mathrm{codim}_{\mathbb{P}^{M-2}} B$ one of
the $q_1, \ldots, q_{j-1}$ must not vanish on $\langle B \rangle$.
Let $i_1$ be the smallest index such that $q_{i_1|\langle B
\rangle} \not\equiv 0$. If $b = j-2$ one of the irreducible
components of $V(q_{i_1}) \cap \langle B \rangle$ must be $B$ and
we are done. For $b < j-2$ choose an irreducible component $R_1$
of $V(q_{i_1}) \cap \langle B \rangle$ containing $B$. Since $\dim
R_1 > \dim B$ and $q_{i|R_1} \equiv 0$ for $i = 1, \ldots, i_1$ we
can find an index $i_1+1 \leqslant i_2 \leq j-1$ such that
$q_{i_2}$ does not vanish on $R_1$: Otherwise, $R_1 \subset V(q_1,
\ldots, q_{j-1})$, and this is a contradiction to
$\mathrm{codim}_{\mathbb{P}^{M-2}} V(q_1, \ldots, q_{j-1}) = j-1$.
In the same way we can inductively find $i_3, \ldots, i_{j-1-b}$
such that finally $B$ is an irreducible component of $V(q_{i_1},
\ldots, q_{i_{j-1-b}}) \cap \langle B \rangle$. Q.E.D. for the
lemma.\vspace{0.1cm}

The lemma implies that $U_j$ is contained in the union of all
locally Zariski-closed subsets $ U_j(P; i_1, \ldots, i_{j-1-b})
\subset \prod_{i=1}^j \mathcal{P}_{m_i}^{M-1}$ of the form
\[ \left\{ (q_1, \ldots, q_j) \left|
   \begin{array}{l}
    \mathrm{there exists an irreducible\ component\ } B\ \mathrm{of\ } V(q_{i_1},
    \ldots, q_{i_{j-1-b}}) \cap P: \\
   \langle B \rangle = P,\ \mathrm{codim}_P B = j-1-b,\  q_{i|B} \equiv 0\
   \mathrm{for\ all\ } i = 1, \ldots, j
   \end{array}
   \right. \right\}, \]
where $P$ ranges over all codimension $b$ linear subspaces of
$\mathbb{P}^{M-2}$, $0 \leqslant b \leqslant j-1$, and the indices
$i_1, \ldots, i_{j-1-b}$ range over all increasing sequences $1
\leqslant i_1 < \cdots < i_{j-1-b} \leqslant j-1$. Note that for
$b = j-1$ we just consider the sets
\[ U_j(P) := \left\{ (q_1, \ldots, q_j) \left|
\begin{array}{l} P\ \mathrm{is\ an\ irreducible\ component}\\
\mathrm{of\ } V(q_1, \ldots, q_{j-1})\ \mathrm{and\ } q_{j|P}
\equiv 0
\end{array} \right. \right\} \subset \prod_{i=1}^j \mathcal{P}_{m_i}^{M-1}. \]
Since the dimension of the Grassmann variety $\mathbb{G}(M-2-b,
M-2)$ is $b(M-1-b)$, estimate (\ref{U_j-deg-2-eq}) will follow if
\begin{equation} \label{UjP-est-eq}
\mathrm{codim}_{\prod_{i=1}^j \mathcal{P}_{m_i}^{M-1}} U_j(P; i_1,
\ldots, i_{j-1-b}) \geqslant 2M + 1 + b(M-1-b).
\end{equation}
holds. By a linear change of coordinates we can assume that
\[ P = \{ z_{M-b} = \cdots = z_{M-1} \}. \]
Let $q_{i_1}, \ldots, q_{i_{j-1-b}}$ be polynomials in
$\prod_{r=1}^{j-1-b} \mathcal{P}_{m_{i_r}}^{M-1}$ such that an
irreducible component $B$ of $V(q_{i_1}, \ldots, q_{i_{j-1-b}})$
lies in $P$ with $\mathrm{codim}_P B = j-1-b$ and $\langle B
\rangle = P$. Then for any degree $m$, products of $m$ linear
forms
\[ \prod_{i=1}^m (a_{i,1}z_1 + \cdots a_{i,M-b-1}z_{M-b-1})\]
do not vanish on $B$. These products span a linear subspace of
$\mathcal{P}_m^{M-1}$ of dimension $(M-b-2)m + 1$ intersecting the
linear subspace of polynomials vanishing on $B$ only in $0$. We
apply these facts to the polynomials $q_i \in
\mathcal{P}_{m_i}^{M-1}$ with $i \in \{1, \ldots, j\} - \{i_1,
\ldots, i_{j-1-b}\}$ and obtain
\begin{eqnarray*}
\mathrm{codim}_{\prod_{i=1}^j \mathcal{P}_{m_i}^{M-1}} U_j(P; i_1,
\ldots, i_{j-1-b}) & \geqslant & \left( \sum_{\substack{i=1 \\ i
\not\in \{i_1, \ldots,
i_{j-b-1}\}}}^j m_i \right)(M-b-2)+ b + 1 \\
 & \geqslant & \left( \sum_{i=1}^b m_i + m_j \right) (M-b-2) + b + 1.
\end{eqnarray*}
So (\ref{UjP-est-eq}) follows from
\begin{equation} \label{m_ib-est-eq}
\left( \sum_{i=1}^b m_i + m_j - b \right) (M-b-2)  \geq 2M.
\end{equation}

Since $m_i \geqslant 2$ for $i=1, \ldots, j-1$ and we consider the
case $m_j \geqslant 3$ we have
\[ \sum_{i=1}^b m_i + m_j - b \geqslant b + 3. \]
The polynomial $(b+3)(M-2-b) - 2M = -b^2 + (M-5)b + M - 6$
quadratic in $b$ increases for $b \leqslant \frac{M-5}{2}$ and
decreases for $b \geqslant \frac{M-5}{2}$. Since $M \geqslant 3k +
4 \geqslant 7$, hence $3(M-2) - 2M \geqslant 0$, estimate
(\ref{m_ib-est-eq}) is shown for $b = 0, M-5$, hence for all $0
\leqslant b \leqslant M-5$. This leaves the cases $b = M-4,
M-3$.\vspace{0.1cm}

If $b = M-4$ then $j = M-3$ or $M-2$, and by the assumptions on
the degrees $m_i$ we have
\[ \sum_{i=1}^b m_i + m_j \geqslant \sum_{i=1}^M m_i - 2d_k - (d_k-1) =
\sum_{l=1}^{k-1} \frac{d_l(d_l+1))}{2} + \frac{(d_k-3)(d_k-2)}{2} + 2 - k. \]
Similarly, if $b = M-3$ then $j = M-2$, and
\[ \sum_{i=1}^b m_i + m_j \geqslant \sum_{i=1}^M m_i - d_k - d_k =
\sum_{l=1}^{k-1} \frac{d_l(d_l+1))}{2} + \frac{(d_k-2)(d_k-1)}{2}
+ 1 - k. \] This last case is the worst possible situation:
$V(q_1, \ldots, q_{M-3})$ is a line in $\mathbb{P}^{M-2}$, and
$q_{M-2}$ vanishes on this line.\vspace{0.1cm}

Solving an optimization problem and using $\sum_{l=1}^k d_l = M+k$
we obtain that
\[ \sum_{l=1}^{k-1} d_l^2 + (d_k-3)^2 \geqslant k
\left( \frac{M-3+k}{k} \right)^2\ \mathrm{and\ }
    \sum_{l=1}^{k-1} d_l^2 + (d_k-2)^2 \geqslant k \left( \frac{M-2+k}{k} \right)^2. \]
Consequently, (\ref{m_ib-est-eq}) follows for $b = M-4, M-3$ if
\begin{equation} \label{M-4-est-eq}
\left[ \frac{(M-3+k)^2}{2k} + \frac{M-3+k}{2} - k - M + 6 \right]
\cdot 2 \geqslant 2M
\end{equation}
and
\begin{equation} \label{M-3-est-eq}
\left[ \frac{(M-2+k)^2}{2k} + \frac{M-2+k}{2} - k - M + 3 \right]
\cdot 1 \geqslant 2M.
\end{equation}

Now (\ref{M-4-est-eq}) is equivalent to
\begin{equation}
k \leqslant \frac{(M-3)^2}{M}
\end{equation}
and (\ref{M-3-est-eq}) is equivalent to
\begin{equation} \label{k-M-3-est-eq}
k \leqslant \frac{(M-2)^2}{3M-2}.
\end{equation}
Both inequalities are satisfied if $M \geqslant 3k+4$. Then $k
\geqslant 1$ implies $M \geqslant 7$, so no further lower bound on
$M$ is needed. Proof of Theorem 7 is complete.\vspace{0.3cm}


{\bf 5. Hypertangent divisors.} Let us prove Theorem 5. The
argument is similar to the proof of Theorem 4 in \cite{Pukh06b},
so we will only sketch the main steps.\vspace{0.1cm}

{\bf Step 1.} Assume that the pair $(V,\frac{1}{n}D)$ is not
canonical for an effective divisor $D\sim nH$, where $H\in
\mathop{\rm Pic} V$ is the class of a hyperplane section,
$K_V=-H$. By linearity of all conditions involved in this
assumption, the divisor $D$ can be assumed to be irreducible. For
some prime divisor $E$ over $V$ the inequality $\mathop{\rm
ord}\nolimits_ED>na(E,V)$ holds.\vspace{0.1cm}

We look at the centre $B\subset V$ of the divisor $E$. Since by
\cite[Proposition 3.6]{Pukh06b} for any irreducible subvariety
$Y\subset V$ of dimension at least $k$ (where $k=\mathop{\rm
codim} (V\subset{\mathbb P})$) we have $\mathop{\rm
mult}\nolimits_Y D\leqslant n$, we conclude that $\mathop{\rm dim}
B\leqslant k-1$.Let $o\in B$ be a point of general position,
$\sigma\colon V^+\to V$ its blow up and $E^+=\sigma^{-1}(o)$ the
exceptional divisor. A standard argument (see, for example,
\cite[Proposition 2.5]{Pukh06b}) gives us a hyperplane
$\Delta\subset E^+\cong {\mathbb P}^{M-1}$ satisfying the
inequality
$$
\mathop{\rm mult}\nolimits_oD+\mathop{\rm
mult}\nolimits_{\Delta}D^+>2n,
$$
where $D^+$ is the strict transform of $D$ on $V^+$.\vspace{0.1cm}

Now let $T$ be a general hyperplane section of $V$, containing the
point $o$ and cutting out $\Delta$ on $E^+$: that is to say,
$T^+\cap E^+=\Delta$. It is easy to see that the effective cycle
of the scheme-theoretic intersection $D_T=(D\circ T)$ is well
defined and satisfies the estimate
\begin{equation}\label{2n}
\mathop{\rm mult}\nolimits_o D_T>2n.
\end{equation}
We will consider $D_T$ as an effective divisor on the complete
intersection $T\subset {\mathbb P}^{M+k-1}$ of codimension $k$,
$D_T\sim nH_T$.\vspace{0.1cm}

{\bf Step 2.} By Sec. 3-4, the complete intersection $T\subset
{\mathbb P}^{M+k-1}$ satisfies the regularity condition. Namely,
for a system of linear coordinates $z_1,\dots, z_{M+k-1}$ with the
origin at the point $o$, the variety $T$ is given by a system of
non-homogeneous polynomial equations:
$$
\begin{array}{lcc}
{\bar f}_1 & = & {\bar q}_{1,1}+\dots + {\bar q}_{1,d_1}, \\
   & \dots &    \\
{\bar f}_i & = & {\bar q}_{i,1}+\dots + {\bar q}_{i,d_i}, \\
   & \dots &    \\
{\bar f}_k & = & {\bar q}_{k,1}+\dots + {\bar q}_{k,d_k},
\end{array}
$$
where the bar means the restriction onto the hyperplane $\{l=0\}$
--- the hyperplane that cuts out $T$ on $V$. Now the set of
homogeneous polynomials
$$
\{{\bar q}_{i,j}\,|\, (i,j)\in I\}
$$
forms a regular sequence in ${\cal O}_{o,{\mathbb
P}^{M+k-1}}$.\vspace{0.1cm}

{\bf Step 3.} Now we can apply the technique of hypertangent
divisors \cite[Chapter 3]{Pukh13a} to the divisor $D_T$ on the
complete intersection $T\subset {\mathbb P}^{M+k-1}$ in precisely
the same way as it was done in \cite[Section 3]{Pukh06b} or, in
more details, in \cite[Section 5]{Pukh10c} and obtain the estimate
$$
\frac{\mathop{\rm mult}\nolimits_o D_T}{2n}\leqslant \mathop{\rm
max}\left\{1,\frac34\cdot\frac{d_k}{d_k-1}\cdot\frac{d_+}{d_+-1}\right\},
$$
where $d_+=d_k$, if $d_{k-1}=d_k$, and $d_+=d_k-1$, otherwise. It
is easy to see that if $d_k\geqslant 8$, this gives us the
inequality $\mathop{\rm mult}\nolimits_o D_T\leqslant 2n$, which
contradicts the estimate (\ref{2n}) obtained in Step 1. The
contradiction completes the proof of Theorem 5.\vspace{0.1cm}

{\bf Proof of Theorem 6}  follows the same lines and repeats the
argument given in \cite[Section 4]{Pukh10c} word for word. What is
different from the proof of Theorem 5 given above, is the starting
point: assuming the inequality
$$
\mathop{\rm ct}(V)\leqslant\frac{M}{M+1},
$$
we obtain for any rational number $\lambda>\frac{M}{M+1}$ an
effective divisor $D\sim nH$ such that the pair
$(V,\frac{\lambda}{n}D)$ is not canonical. Now, repeating the
proof of Theorem 5, we use the technique of hypertangent divisors
to obtain the inequality
$$
1<\lambda\mathop{\rm
max}\left\{1,\frac34\cdot\frac{d_k}{d_k-1}\cdot\frac{d_+}{d_+-1}\right\},
$$
and taking the limit as $\lambda\to \frac{M}{M+1}$, we conclude
that
$$
\mathop{\rm
max}\left\{1,\frac34\cdot\frac{d_k}{d_k-1}\cdot\frac{d_+}{d_+-1}\right\}
\geqslant\frac{M+1}{M}.
$$
However, it is a trivial check that in the assumptions of Theorem
6 the last inequality can not be true. Q.E.D. for Theorem 6.

\begin{flushleft}
Department of Mathematical Sciences\\
The University of Liverpool\\
England
\end{flushleft}

\noindent  e-mail: eckl@liv.ac.uk,\quad pukh@liv.ac.uk

\end{document}